\documentclass[11pt]{article}
\usepackage{amsmath,amsfonts,amsbsy,amsgen,amscd, amsxtra,amstext,amsopn,amssymb}
\usepackage{latexsym}
\usepackage{theorem}
\usepackage{mathrsfs}
\usepackage{pictex}
\usepackage{epic}
\usepackage{euscript}

{\theorembodyfont{\slshape}
\newtheorem{theorem}{Theorem}[section]
\newtheorem{proposition}[theorem]{Proposition}
\newtheorem{lemma}[theorem]{Lemma}
\newtheorem{corollary}[theorem]{Corollary}
}
{\theorembodyfont{\rmfamily}

\newtheorem{remark}[theorem]{Remark}

}

\def\proof{{\noindent\sc Proof. \quad}}
\newcommand\proofof[1]{{\medskip\noindent\sc Proof of #1.\quad}}
\def\eproof{{\mbox{}\hfill\qed}\medskip}
\newcommand\qed{{\unskip\nobreak\hfil\penalty50\hskip2em\vadjust{}
\nobreak\hfil$\Box$\parfillskip=0pt\finalhyphendemerits=0\par}}

\def\R{\mathbb{R}}

\def\E{\mathop\mathbb{E}}

\def\Avr{\mathbf{Avr}}

\newcommand{\e}{\varepsilon}

\newcommand{\s}{\sigma}

\newcommand{\z}{\zeta}

\def\Prob{\mathop{\rm Prob}}
\newcommand{\oA}{\overline{A}}
\newcommand{\Id}{\mathrm{I}}
\newcommand{\Sph}{\mathbb{S}}
\newcommand{\spann}{\mathsf{span}}
\newcommand{\vol}{\mathsf{vol\,}}
\newcommand{\kap}{\mathsf{cap}}

\newcommand{\LoP}{\mathsf{LoP}}
\newcommand{\transp}{^{\rm T}}
\def\Oh{{\cal O}}

\newcommand{\bOA}{{{\mathbf{Avr}\left(\ln\kappa(A)\right)}}}
\newcommand{\bNA}{{\mu(m,n,\sigma)}}

\newcommand{\size}{\mathsf{size}}

\begin{document}

\begin{title}
{\Large {\bf Smoothed Analysis of Moore-Penrose Inversion}}
\end{title}
\author{Peter B\"urgisser\thanks{Institute of Mathematics, University of
Paderborn, Germany. Partially supported
by DFG grant BU 1371/2-1 and BU 1371/3-1.}\ \ and
Felipe Cucker\thanks{Dept.\ of
Mathematics, City University of Hong Kong,
Kowloon Tong, Hong Kong. Partially supported by a grant
from the Research Grants Council of Hong Kong,
project No.\ CityU 100808.}}
\date{}
\makeatletter
\maketitle

\begin{quote}{\small
{\bf Abstract.}
We perform a smoothed analysis of the condition number
of rectangular matrices. We prove that, asymptotically, the
expected value of this condition number depends only of the
elongation of the matrix, and not on the center and variance
of the underlying probability distribution.
}\end{quote}


\section{Introduction}\label{se:intro}

The most widely used extension to rectangular matrices
of the notion of inverse of square matrices is
the so called Moore-Penrose inverse. For a full rank
matrix $A\in\R^{m\times n}$ this is defined as
$A^\dagger:=(A\transp A)^{-1}A\transp$ if $m\geq n$,
and as $A^\dagger:=A\transp (AA\transp)^{-1}$, otherwise.
Immediate applications of $A^\dagger$ include the solution
of least square problems
\begin{equation}\label{eq:ls}
    \min_{x\in\R^n}\|Ax-b\|^2,
\end{equation}
with $b\in\R^m$ and $m>n$, or of smallest solutions of
underdetermined systems
\begin{equation}\label{eq:us}
    \min_{x\mid Ax=b}\|x\|^2
\end{equation}
when $n>m$.
In both cases, the solution is given by $x=A^\dagger b$.
Well known results in error analysis show that the accuracy in the
computation of $A^\dagger$, or in the computation of the solution
$x$ for the problems above, crucially depends on the {\em condition
number} $\kappa(A):=\|A\|\,\|A^\dagger\|$ of $A$, where $\|A\|$
denotes the spectral norm (see~\cite[Ch.~19]{higham:96}). Accuracy
analysis is not the only source of interest in $\kappa(A)$.
Algorithms such as the conjugate gradient method produce approximate
solutions of linear systems $Px=c$ ---here $P\in\R^{m\times m}$ is a
positive definite matrix and $c\in\R^m$--- with a number of
iterations proportional to $\sqrt{\kappa(P)}$ and, in many cases,
the matrix $P$ has been obtained as $P=AA\transp$ for some matrix
$A\in\R^{m\times n}$. In those cases, $\sqrt{\kappa(P)}=\kappa(A)$
and one is again interested in the latter, this time by complexity
considerations.

The condition number $\kappa(A)$ is not directly readable
from $A$, and its computation seems to require that
of $A^\dagger$. This is a common situation in numerical
analysis. A way out of it, proposed as early as 1951
by von Neumann and Goldstine~\cite{vNGo51} and
more recently pioneered by
Demmel~\cite{Demmel88} and Smale~\cite{Smale97},
consists of randomizing the matrix~$A$ ---say, by endowing
$\R^{m\times n}$ with a multivariate standard Gaussian
distribution $N(0,\Id)$--- and considering
its condition number as a derived random variable.

In Chen and Dongarra~\cite{ChDo:05} the
following tail estimates on $\kappa(A)$ were shown
for $A\in\R^{m\times n}$ with $n\ge m$:
for $x\ge n-m+1$ we have
\begin{equation}\label{eq:chen-dong}
 \frac1{\sqrt{2\pi}}\, \Big(\frac{1}{5x}\Big)^{n-m+1}\ \leq\
 \Prob_{A\sim N(0,\Id)}
 \Big\{\kappa(A)\geq\frac{x}{1-\lambda}\Big\}\
  \leq\ \frac1{\sqrt{2\pi}}\, \Big(\frac{7}{x}\Big)^{n-m+1} .
\end{equation}
Moreover, the expectation $\E(\kappa(A))$ can be bounded as a function
of the {\em elongation} $\frac{m-1}{n}$ only, independently of~$n$.
(We remark that this is not true for Demmel's scaled
condition number $\|A\|_F\,\|A^\dagger\|$, compare~\cite{edelm:92}.)
More precisely,
for a sequence $(m_n)$ of integers such that
$\lim_{n\to\infty}m_n/n=\lambda\in(0,1)$ and a
sequence of standard Gaussian random matrices $A_n \in\R^{m_n\times n}$,
we have in almost sure convergence
\begin{equation}\label{eq:Edelman_kappa}
  \kappa(A_n)\stackrel{\rm a.s.}{\longrightarrow}
  \frac{1+\sqrt{\lambda}}{1-\sqrt{\lambda}} .
\end{equation}
This follows from Geman~\cite{geman:80} and Silverstein~\cite{silver:85}
(see Edelman~\cite{edelm:88} for more precise results).

The above results provide theoretical reasons of why least squares
problems such
as~\eqref{eq:ls} or underdetermined systems such as~\eqref{eq:us}
are solved to great accuracy or why the conjugate gradient method is
so efficient in practice. In fact, it follows from~\eqref{eq:Edelman_kappa}
that the expected number of iterations of
the conjugate gradient method on the random input $P=AA\transp $ remains
bounded in terms of the {\em elongation} $m/n$ as $n\to\infty$ and
$A\in\R^{m\times n}$ is standard Gaussian. Our main result stated
below implies that this phenomenon is still true for any matrix
that is only slightly perturbed.

The choice of $N(0,\Id)$ as underlying data distribution
is pervasive in the {\em average-case analysis} of condition
numbers
(and other quantities occurring in numerical analysis). It has the
virtue of simplicity as a first approach to understanding
which condition numbers one may expect. But it has been
criticized due to the loose relationship of the Gaussian
$N(0,\Id)$ to the measures that may be governing
data drawing in practice.
In particular, it has been observed
that the use of Gaussians may be `optimistic' in the sense
that they may put more probability mass on the instances where
the values of the function $\psi$ under consideration are small.
Such an optimism would produce yield an expectation
$\E(\psi)$ smaller than the true one.

An alternate, more conservative,
form of analysis has been proposed by Spielman and Teng
under the name of {\em smoothed analysis}. It replaces
the Gaussian measure $N(0,\Id)$ by the measures
$N(\oA,\sigma^2\Id)$ where $\oA$ is arbitrary. The idea is
then to replace the unlikely `average data' by a (usually
small) perturbation of any possible occurring data. The
rationale for this form of analysis is offered in a number
of papers~\cite{ST-simplex:04,sst:06,ST:06,ST:09} and
we won't repeat it here in full.
We note, nonetheless,
that the local nature of randomization in smoothed analysis,
coupled with its worst-case dependence on the input data,
removes from smoothed analysis the possible optimism
we mentioned above for average-case analysis.
In recent years,
different aspects of algorithm behavior for a variety of
problems have been analyzed this way. These include
condition numbers of square matrices with
real~\cite{Wsch:04} or $\{-1,1\}$
coefficients~\cite{TaoVu:07}, complexity of
interior-point methods~\cite{DuSpTe:09}, and machine
learning~\cite{ArVa:09}. The typical satisfying result
is {\em polynomial smoothed complexity}
(see~\cite[Def.~2]{ST:09}), consisting of a bound
of the form
\begin{equation}\label{eq:sc}
  \sup_{\oA} \E_{A\sim N(\oA,\sigma^2\Id)} \psi(A)
  \leq c \sigma^{-k_1} \size(A)^{k_2}
\end{equation}
where $\psi$ is the function whose behavior we are
analyzing and $c,k_1,k_2$ are positive constants.

In this paper we provide a smoothed analysis for
Moore-Penrose inversion, extending~\eqref{eq:chen-dong}
from the average-case analysis to smoothed analysis.
To state the results we need to introduce some notations.
We assume $1\le m\le n$ throughout the paper.
For a standard Gaussian $X\in\R^{m\times n}$ we put
\begin{equation}\label{eq:defQ}
 Q(m,n) := \frac1{\sqrt{n}}\, \E(\|X\|) .
\end{equation}
(Lemma~\ref{pro:Espnorm} shows that $Q(m,n) \le 6$ .)
We define for $\lambda \in (0,1)$ the quantity
\begin{equation}\label{eq:defc}
 c(\lambda) \ :=\   \sqrt{\frac{1+\lambda}{2(1-\lambda)}} .
\end{equation}
Note that $c(\lambda)$ is monotonically increasing,
$\lim_{\lambda\to 0}c(\lambda) = \frac1{\sqrt{2}}$ and
$\lim_{\lambda\to 1}c(\lambda) = \infty$.
Further, for $1\le m\le n$ and $0<\sigma\le 1$, we
define the {\em elongation}
$\lambda := \frac{m-1}{n}$ and introduce the quantity
\begin{equation}\label{eq:defB}
 \z_\s(m,n) \ :=\ \Big(Q(m,n) + \frac1{\s\sqrt{n}}\Big)\
  c(\lambda)^{\frac1{n-m+1}}  .
\end{equation}

Our main result is the following tail bound on the condition
number of rectangular matrices under local Gaussian perturbations.

\begin{theorem}\label{thm:tailbound}
Suppose that $\oA\in\R^{m\times n}$ satisfies
$\|\oA\| \le 1$ and let $0<\sigma\le 1$.
Put $\lambda := \frac{m-1}{n}$.
Then, for $z \ge \z_\s(m,n)$, we have 
$$
  \Prob_{A\sim N(\oA,\sigma^2\Id)}
   \Big\{\kappa(A)\geq\frac{ez}{1-\lambda}\Big\}\
  \leq\
   2c(\lambda)\bigg[
   \Big(Q(m,n) + \sqrt{2\ln(2z)} + \frac1{\s\sqrt{n}}\Big)\,
    \frac1z\bigg]^{n-m+1} .
$$
\end{theorem}

\begin{remark}
{\bf 1.} The decay in~$z$ in this tail bound is the same as
in~\eqref{eq:chen-dong} up to the logarithmic factor $\sqrt{\ln z}$.
We believe that the latter is an artefact of our proof that could be omitted.
In fact, the exponent $n-m+1$ is just the codimension of the set
$\Sigma := \{A\in\R^{m\times n} \mid \mathrm{rk} A < m\}$
of rank deficient matrices, cf.~\cite{harr:95}.
Moreover, it is known~\cite{higham:96} that
$\|A^\dagger\| = 1/\mathrm{dist}(A,\Sigma)$
where the distance is measured in the Euclidean norm.
From the interpretation of $\Prob\{\kappa(A)\geq t\}$
as the volume of a tube around $\Sigma$, as discussed in~\cite{BCL:06a},
one would therefore expect a decay of order $1/z^{n-m+1}$.

{\bf 2.} When $\sigma=1$ and $\oA=0$, Theorem~\ref{thm:tailbound}
yields tail bounds for the usual average case. One may therefore
compare these bounds with~\eqref{eq:chen-dong}.
In doing so, we see that
the bound in Theorem~\ref{thm:tailbound}
has the additional factor~$c(\lambda)$ (going to $\infty$
for $\lambda\to 1$). However, we note that the bound~\eqref{eq:chen-dong}
holds only for
$x=ez\ge n-m + 1$,
while our bound holds for any $z \ge \z_\s(m,n)$.
Furthermore, if we fix
$\lambda\in(0,1)$ and let $(m_n)$ be a sequence of positive integers
such that $\lim m_n/n=\lambda$,
it follows from~\cite{geman:80} that
$$
\lim_{n\to\infty} Q(m_n,n)= 1 + \sqrt{\lambda}.
$$
This implies that
$\lim_{n\to\infty} \z_\s(m_n,n) = 1 + \sqrt{\lambda}$
for fixed $\s \in (0,1]$ and, in particular, that
$\z_\s(m_n,n) \le 2$ for sufficiently large~$n$ .
That is, for large $n$, the tail bound in Theorem~\ref{thm:tailbound}
is valid for any $z\ge 2$.
\end{remark}

Theorem~\ref{thm:tailbound} easily implies the following
bound on expectations.

\begin{corollary}\label{cor:main}
For all $\lambda_0\in(0,1)$ 
there exists $n_0$ such that for all
$1\le m\le n$ such that
$\lambda = \frac{m-1}{n} \le \lambda_0$ and
$n\ge n_0$ we have
for all $\sigma$ with $\frac1{\sqrt{m}}\le \sigma\le 1$,
and all
$\oA\in\R^{m\times n}$ with $\|\oA\| \le 1$, that
$$
 \E_{A\sim N(\oA,\sigma^2\Id)}(\kappa(A))
  \ \leq\ \frac{20.1}{1-\lambda} .
$$
\end{corollary}

As for the average-case analysis, this bound
is independent of~$n$ and depends only on the bound $\lambda_0$
on the elongation.
Thus we have a bound of type~\eqref{eq:sc} with $k_2=0$.
Surprisingly, the smoothed complexity bound in Corollary~\ref{cor:main}
is also independent of~$\sigma$. We thus add reasons
---and we will become more specific in
Section~\ref{sec:applications}---
to the current understanding of the accuracy
in least squares or underdetermined system solving
or the complexity of the conjugate gradient method.

A first approach to the smoothed analysis of Moore-Penrose
inversion appears in~\cite{CDW:05}. The bounds obtained in
that paper are worse by an order of magnitude
than those we obtain here.
In Section~\ref{sec:numerical} we compare these bounds
{\sf with ours} 
as well as with actual averages obtained,
for specific values of $n,m$ and $\sigma$, in numerical
simulations.

Our proof techniques are an extension of methods
employed by Sankar et al.~\cite{sst:06}.

\medskip

\noindent{\bf Acknowledgements.}
This work was carried out during the special semester on Foundations
of Computational Mathematics in the fall of 2009.
We thank the Fields Institute in Toronto for hospitality and
financial support.

\section{Preliminaries}

\subsection{Some definitions and notation}

The {\em spectral norm} of a matrix $A\in\R^{m\times n}$ is defined as
$\|A\| := \sup_{\|x\|=1} \|Ax\|$, where $\|x\|$ denotes the Euclidean norm.
The {\em Frobenius norm} of $A$ is defined as the Euclidean norm of $A$
when interpreted as a vector.

Suppose that $A\in\R^{m\times n}$ is of maximal rank and $m\le n$.
The {\em Moore-Penrose inverse} of $A$
is defined as
$A^\dagger := A\transp  (AA\transp)^{-1}\in\R^{n\times m}$.
It can also be characterized as follows.
For any $v\in\R^m$ the vector $w=A^\dagger v$ is orthogonal to
the kernel of~$A$ and satisfies $Aw=v$.
The {\em condition number} $\kappa(A)$ is defined as
$\kappa(A) :=\|A\| \cdot \|A^\dagger\|$.

Let $\oA\in \R^{m\times n}$ and $\s>0$.
The {\em isotropic normal distribution} $N(\oA,\s \Id)$ with center~$\oA$
and covariance matrix $\s^2\Id$ is the probability distribution on
$\R^{m\times n}$ with the density
$$
 \rho_{\oA,\s}(A) := \frac1{(2\pi)^{\frac{mn}{2}}}\,
  e^{-\frac{\|A-\oA\|_F^2}{2\s^2}} .
$$

\begin{lemma}\label{le:gbound}
For $\lambda\in (0,1)$ we have
$\lambda^{-\frac{\lambda}{1-\lambda}} \le e$.
\end{lemma}

\proof
Writing $u=1/\lambda$ the assertion is equivalent to
$u^{\frac1{u-1}} \le e$ or $u\le e^{u-1}$,
which is certainly true for $u \ge 1$.
\eproof

\subsection{Concentration on spheres}

Let $\Sph^{m-1}:=\{x\in\R^m \mid \|x\|=1\}$ denote the unit sphere in $\R^m$.
We denote by $\Oh_{m-1}$ its volume, which is given by
$\Oh_{m-1} = 2\pi^{m/2}/\Gamma(\frac{m}{2})$.

The following estimate tells us how likely a random point on $\Sph^{m-1}$
will lie in a fixed spherical cap.

\begin{lemma}\label{lem:s}
Let $u\in\Sph^{m-1}$ be fixed, $m\geq 2$. Then, for all $\xi\in[0,1]$,
$$
   \Prob_{v\sim U(\Sph^{m-1})}
   \big\{\big|u\transp v\big|\geq \xi\big\} \geq
   \sqrt{\frac{2}{\pi m}}\; (1-\xi^2)^{\frac{m-1}{2}}.
$$
\end{lemma}

\proof
We put $\theta=\arccos\xi$ and let $\kap(u,\theta)$ denote
the spherical cap in $\Sph^{m-1}$ with center~$u$ and angular radius~$\theta$.
Using the bounds in Lemmas~2.1 and 2.2 of~\cite{bcl:08a} we get
$$
   \Prob_{v\sim U(\Sph^{m-1})}
   \big\{\big|u\transp v\big|\geq \xi\big\}
  = \frac{2\,\vol\kap(u,\theta)}{\vol\Sph^{m-1}}
  \geq \frac{2\Oh_{m-2}}{\Oh_{m-1}}\;
   \frac{(1-\xi^2)^{\frac{m-1}{2}}}{(m-1)}.
$$
Using the formula for $\Oh_{m-1}$ and the recursion
$\Gamma(x+1)=x\Gamma(x)$ we have
$$
    \frac{\Oh_{m-2}}{\Oh_{m-1}}
  = \frac{1}{\sqrt{\pi}}
   \frac{\Gamma\big(\frac{m}{2}\big)}
   {\Gamma\big(\frac{m-1}{2}\big)}
  = \frac{1}{\sqrt{\pi}}
   \frac{\Gamma\big(\frac{m+1}{2}\big)}
   {\Gamma\big(\frac{m-1}{2}\big)}
   \frac{\Gamma\big(\frac{m}{2}\big)}
   {\Gamma\big(\frac{m+1}{2}\big)}
   = \frac{m-1}{2\sqrt{\pi}}
   \frac{\Gamma\big(\frac{m}{2}\big)}{\Gamma\big(\frac{m+1}{2}\big)}.
$$
The assertion follows now from the estimate
\begin{equation}\label{eq:Gine}
 \frac{\Gamma\big(\frac{m}{2}\big)}{\Gamma\big(\frac{m+1}{2}\big)}
 \ \ge\ \sqrt{\frac{2}{m}}.
\end{equation}
This estimate can be quickly seen as follows.
Suppose that $Z\in\R^{m}$ is standard normal distributed.
Using polar coordinates and the variable
transformation $u=\rho^2/2$ we get
\begin{eqnarray}\notag
 \E(\|Z\|) &= &\frac{\Oh_{m-1}}{(2\pi)^{\frac{m}{2}}}
 \int_0^\infty \rho^m e^{-\frac{\rho^2}{2}} d\rho
 \:=\: \frac{\Oh_{m-1}}{(2\pi)^{\frac{m}{2}}}\, 2^{\frac{m-1}{2}}
    \int_0^\infty u^{\frac{m-1}{2}} e^{-u} du \\ \label{eq:remember}
 &=& \frac{\Oh_{m-1}}{(2\pi)^{\frac{m}{2}}}\, 2^{\frac{m-1}{2}}
    \Gamma(\frac{m+1}{2})
 \:=\: \sqrt{2}\,\frac{\Gamma(\frac{m+1}{2})}
   {\Gamma(\frac{m}{2})},
\end{eqnarray}
where we used the definition of the Gamma function
for the second last equality.
To complete the proof of \eqref{eq:Gine} we note that
$\E(\|Z\|) \le \sqrt{\E(\|Z\|^2)} =\sqrt{m}$.
\eproof

For later use we note that~\eqref{eq:remember} implies
$$
 \frac{\Gamma\big(\frac{m+1}{2}\big)}{\Gamma\big(\frac{m}{2}\big)}
 = \frac{\Gamma\big(\frac{m+2}{2}\big)}{\Gamma\big(\frac{m}{2}\big)}
      \frac{\Gamma\big(\frac{m+1}{2}\big)}{\Gamma\big(\frac{m+2}{2}\big)}
 = \frac{m}{2}\, \frac{\Gamma\big(\frac{m+1}{2}\big)}{\Gamma\big(\frac{m+2}{2}\big)}
 \ \ge\  \frac{m}{2}\, \sqrt{\frac{2}{m+1}},
$$
using~\eqref{eq:Gine} for the right-hand inequality. Therefore
\begin{equation}\label{eq:chi-lb}
 \E(\|Z\|)\ \ge\ \frac{m}{\sqrt{m+1}} .
\end{equation}

\subsection{Large deviations}

We will use a powerful large deviation result.
Let $F:\R^N\to\R$ be a Lipschitz continous function with Lipschitz
constant~$L$, so that $|F(x)-F(y)| \le L \|x-y\|$ for all $x,y\in\R^N$,
where $\|\ \|$ denotes the Euclidean norm.
Now suppose that $x\in\R^N$ is a standard Gaussian random vector such that
$\E(F(x))$ exists. Then it is known~\cite[(1.4)]{letala:91} that for all $t>0$
\begin{equation}\label{eq:concentration}
  \Prob\{F(x)\geq \E(F) +t\}\leq e^{-\frac{t^2}{2L^2}}.
\end{equation}
(We note that in ~\cite[(1.4)]{letala:91} this is only stated for the median, but
the inequality holds as well for the expectation. See also~\cite{ledoux:01}.)

\subsection{A bound on the expected spectral norm}

The function $\R^{m\times n}\to\R$ mapping a matrix $X$ to its
spectral norm~$\|X\|$ is Lipschitz continuous with Lipschitz constant~$1$,
as $\|X-Y\| \le \|X-Y\|_F$.
The concentration bound~\eqref{eq:concentration}, together
with~\eqref{eq:defQ}, implies that for $t>0$,
\begin{equation}\label{eq:spnb}
\Prob\Big\{\|X\| \ge Q(m,n)\sqrt{n} + t\}\ \le\ e^{-\frac{t^2}{2}} .
\end{equation}
This tail bound easily implies the following large deviation result.

\begin{proposition}\label{le:enorm}
Let $\oA\in\R^{m\times n}$ with $m\le n$, $\|\oA\| \le 1$, and $\sigma\in(0,1]$.
If $A\in\R^{m\times n}$ follows the law $N(\oA,\sigma^2\Id)$, then, for $t>0$,
$$
 \Prob_{A\sim N(\oA,\sigma^2\Id)}\Big\{\|A\|\geq Q(m,n)\sigma\sqrt{n}+t+1\Big\}
 \ \leq\  e^{-\frac{t^2}{2\sigma^2}}.
$$
\end{proposition}

\proof
We note that
$\|A\|\geq Q(m,n)\sigma\sqrt{n}+t+1$ implies that
$\|A -\oA\| \geq \|A\| -\|\oA\| \ge Q(m,n)\sqrt{n} + t$.
Moreover, if $A\in\R^{m\times n}$ follows the law $N(\oA,\sigma^2\Id)$, then
$X:=\frac{A-\oA}{\sigma}$ is standard Gaussian in $\R^{m\times n}$.
The assertion follows from~\eqref{eq:spnb}.
\eproof

We derive now an upper bound on $Q(m,n)$.
Such result should be well-known but we could not locate in the literature.

\begin{lemma}\label{pro:Espnorm}
For $n>1$ we have
$\sqrt{\frac{n}{n+1}}
\le Q(m,n)\le
 2\Big( 1 + \sqrt{\frac{2\ln (2m-1)}{n}} + \frac1{\sqrt{n}}\Big) \le 6$.
\end{lemma}

The proof relies on the following lemma.

\begin{lemma}\label{le:maxchi}
Let $r_1,\ldots,r_n$ be independent random variables with nonnegative values
such that $r_i^2$ is $\chi^2$-distributed with $f_i$ degrees of freedom.
Then,
$$
 \E\Big(\max_{1\le i \le n} r_i\Big) \le \max_{1\le i\le n}\sqrt{f_i}
   + \sqrt{2\ln n} +  1 .
$$
\end{lemma}

\proof
We start by a large deviation estimate for $\chi^2$-distributed random variables.
Note that $\R^f\to\R$, $x\mapsto \|x\|$, is Lipschitz continuous with
Lip\-schitz constant~$1$.
From \eqref{eq:concentration} we know that for standard Gaussian $x\in\R^n$
and all $t>0$,
$$
  \Prob\{\|x\| \geq \E(\|x\|) +t\}\ \leq\ e^{-\frac{t^2}{2}}.
$$
Since $\E(\|x\|)\le \sqrt{\E(\|x\|^2)} = \sqrt{f}$, this implies for all $t>0$,
\begin{equation}\label{eq:chi-conc}
  \Prob\{\|x\| \geq \sqrt{f} +t \}\ \leq\ e^{-\frac{t^2}{2}}.
\end{equation}

We suppose now that $r_1,\ldots,r_n$ are
independent random variables with nonnegative values
such that $r_i^2$ is $\chi^2$-distributed with $f_i$ degrees of freedom.
Put $f:=\max_i f_i$.
Equation~\eqref{eq:chi-conc} tells us that for all $i$ and all $t>0$,
$$
 \Prob\{ r_i \ge \sqrt{f} + t \}\ \le\  e^{-\frac{t^2}{2}}
$$
and hence, by the union bound,
$$
 \Prob\Big\{\max_{1\le i\le n} r_i \ge \sqrt{f} + t \Big\}\
  \le\ n e^{-\frac{t^2}{2}} .
$$
For a fixed parameter~$b\ge 1$ (to be determined later), this implies
\begin{eqnarray*}
 \E(\max_{1\le i\le n} r_i) &\le& \sqrt{f} + b + \int_{\sqrt{f} + b}^\infty
                       \Prob\{\max_{1\le i\le n} r_i \ge T \}\, dT \\
  &=& \sqrt{f} + b + \int_{b}^\infty \Prob\{\max_{1\le i\le n}  r_i \ge \sqrt{f} + t \}\, dt \\
  &\le& \sqrt{f} + b + n \int_{b}^\infty  e^{-\frac{t^2}{2}}\,dt .
\end{eqnarray*}
Using the well-known estimate
$$
 \frac1{\sqrt{2\pi}}\int_b^\infty  e^{-\frac{t^2}{2}}\,dt
 \le \frac1{b\sqrt{2\pi}}\,e^{-\frac{b^2}{2}}
  \le \frac1{\sqrt{2\pi}}\,e^{-\frac{b^2}{2}}
$$
we obtain
$$
 \E(\max_{1\le i\le n} r_i) \le \sqrt{f} + b + n e^{-\frac{b^2}{2}} .
$$
Finally, choosing $b:=\sqrt{2\ln n}$ we get
$$
 \E(\max_{1\le i\le n} r_i) \le \sqrt{f} + \sqrt{2\ln n} + 1,
$$
as claimed.
\eproof

\proofof{Lemma~\ref{pro:Espnorm}}
A general matrix $X\in\R^{m\times n}$ can be transformed into a
bidiagonal matrix of the form
$$
Y:=\begin{bmatrix}
v_n    &       &         &      & 0     &\cdots &0     \\
w_{m-1}&v_{n-1}&         &      & \vdots&       &\vdots\\
       &\ddots &\ddots   &      & \vdots&       &\vdots\\
       &       & w_1     &v_{n-m+1}&0   &\cdots &0
\end{bmatrix}
$$
with $v_i,w_j\ge 0$ by performing Householder transformations
from the left and right hand side of~$X$, cf.~\cite[\S5.4.3]{golloan:83}.
In particular, $\|X\| = \|Y\|$.
An analysis of this transformation shows that if we start
with a standard Gaussian matrix $X$, then
the $v_n,\ldots,v_{n-m+1},w_{m-1},\ldots,w_1$
are independent random variables such that $v_i^2$ and $w_i^2$ are
$\chi^2$-distributed with $i$~degrees of freedom, cf.~\cite{silver:85}.

The spectral norm of $Y$ is bounded by
$\max_{i}v_i + \max_j w_j \le 2 r$, where
$r$~denotes the maximum of the values $v_i$ and $w_j$.
Lemma~\ref{le:maxchi} implies that, for $n>1$,
$$
 \E(r) \le \sqrt{n} + \sqrt{2\ln (2m-1)} + 1
  \le 3\sqrt{n} .
$$
This shows the claimed upper bound on $Q(m,n)$.
For the lower bound we note that
$\|Y\| \ge |v_n|$ which gives
$\E(\|Y\|) \ge \E(|v_n|)$.
The claimed lower bound now follows from~\eqref{eq:chi-lb},
which states that
$\E(|v_n|) \ge \sqrt{\frac{n}{n+1}}$.
\eproof

\section{Proof of the main results}

The main work consists of deriving tail bounds on $\|A^\dagger\|$,
which is done in the next subsection.

\subsection{Tail bounds for $\|A^\dagger\|$}\label{se:tb-Adag}

\begin{proposition}\label{thm:main_tail}
Let $\oA\in\R^{m\times n}$, $\s>0$, and put
$\lambda :=\frac{m-1}{n}$.
For random $A\sim N(\oA,\s^2\Id)$
we have, for any $t>0$,
$$
 \Prob_{A\sim N(\oA,\s^2\Id)}\Big\{\|A^\dagger\|\geq \frac{t}{1-\lambda}\Big\}
  \ \leq\ c(\lambda)\,
       \bigg(\frac{e}{\s\sqrt{n}\, t}\bigg)^{(1-\lambda)n}.
$$
\end{proposition}

We first show the following result.

\begin{proposition}\label{prop:bound1}
For all $v\in \Sph^{m-1}$,
$\oA\in\R^{m\times n}$, $\s>0$, and $\xi>0$
we have
$$
  \Prob_{A\sim N(\oA,\s^2\Id)}
  \big\{\|A^\dagger v\|\geq \xi \big\}\ \leq\
       \frac{1}{(\sqrt{2\pi})^{n-m+1}}\,
               \frac{\Oh_{n-m}}{n-m+1}\,
               \Big(\frac{1}{\s \xi}\Big)^{n-m+1}.
$$
\end{proposition}

\proof
We first claim that, because of unitary invariance,
we may assume that $v=e_m:=(0,\ldots,0,1)$.
To see this, take $\Phi\in U(m)$ such that $v=\Phi e_m$.
Consider the isometric map
$A\mapsto B=\Phi^{-1}A$ which transforms the density
$\rho_{\oA,\s}(A)$ into a density of the same form, namely
$\rho_{\Phi^{-1}\oA,\s}(B)$.
Thus the assertion for $e_m$ and random $B$
implies the assertion for $v$ and $A$, noting that
$A^\dagger v=B^\dagger e_m$. This proves the claim.

We are going to characterize the norm of $w:=A^\dagger e_m$
in a geometric way.
Let $a_i$ denote the $i$th row of $A$. Almost surely,
the rows $a_1,\ldots,a_m$ are linearly independent;
hence, we assume so in what follows.
Let
$$
  R:=\spann\{a_1,\ldots,a_{m}\},\ S:=\spann\{a_1,\ldots,a_{m-1}\} .
$$
Let $S^\perp$ denote the orthogonal
complement of $S$ in $\R^n$. We decompose
$a_m=a_m^\perp + a_m^S$,
where $a_m^\perp$ denotes the orthogonal projection
of $a_m$ onto $S^\perp$ and $a_m^S\in S$. Then $a_m^\perp\in R$
since both $a_m$ and $a_m^S$ are in $R$. It follows that
$a_m^\perp\in R\cap S^\perp$.

We claim that $w\in R\cap S^\perp$ as well. Indeed,
note that $R$ equals the orthogonal complement of the
kernel of $A$ in $\R^n$. Therefore, by definition of the
Moore-Penrose inverse, $w=A^\dagger e_m$ lies in $R$.
Moreover, since $A A^\dagger=\Id$, we have
$\langle w,a_i\rangle =0$ for $i=1,\ldots,m-1$
and hence $w\in S^\perp$ as well.

It is immediate to see that $\dim R\cap S^\perp=1$.
It then follows that $R\cap S^\perp=\R w = \R a_m^\perp$.
Since $\langle w,a_m\rangle =1$, we get
$1 = \langle w,a_m\rangle = \langle w,a_m^\perp \rangle
      = \|w\|\, \|a_m^\perp\|$
and therefore
\begin{equation}\label{eq:star}
 \|A^{\dagger}e_m\| = \frac1{\|a_m^\perp\|}.
\end{equation}

Let $A_m\in\R^{(m-1)\times n}$ denote the matrix
obtained from~$A$ by omitting $a_m$. The density
$\rho_{\oA,\s}$ factors as
$\rho_{\oA,\s}(A) =\rho_1(A_n)\rho_2(a_n)$ where
$\rho_1$ and $\rho_2$ denote the density functions of
$N(\oA_m,\s^2\Id)$ and
$N(\bar{a}_m,\s^2\Id)$, respectively
(the meaning of $\oA_m$ and $\bar{a}_m$
being clear). Fubini's Theorem combined
with \eqref{eq:star} yield, for $\xi>0$,
\begin{eqnarray}\label{eq:WX}
  \Prob_{N(\oA,\s^2\Id)}\big\{\|A^{\dagger}e_m\|\geq \xi\big\}
   &=&
     \int_{\|A^{\dagger}e_m\| \ge \xi}
     \rho_{\oA,\s^2\Id}(A)\, dA\\
  &=&
     \int_{A_{m}\in\R^{(m-1)\times n}}
     \rho_1(A_{m})
     \cdot  \left(\int_{\|a_m^\perp\|\leq 1/\xi}
     \rho_2(a_m)\, da_m\right) dA_{m}.\nonumber
\end{eqnarray}
To complete the proof it is sufficient to show the bound
\begin{equation}\label{eq:C4}
    \int_{\|a_m^\perp\|\leq\frac1\xi}
    \rho_2(a_m)\, da_m
    \leq   \frac{1}{(\sqrt{2\pi})^{n-m+1}} \,
               \frac{\Oh_{n-m}}{n-m+1}\,
               \Big(\frac{1}{\s \xi}\Big)^{n-m+1}
\end{equation}
for fixed, linearly independent $a_1,\ldots,a_{m-1}$ and $\xi>0$.

To show~\eqref{eq:C4} note that
$a_m^\perp\sim N(\bar{a}_m^\perp,\s^2\Id)$
in $S^\perp\simeq\R^{n-m+1}$
where $\bar{a}_m^\perp$
is the orthogonal projection of $\bar{a}_m$ onto
$S^\perp$.
Let $B_r$ denote the ball of radius~$r$ in $\R^p$ centered
at the origin. It is easy to see that $\vol B_r = \Oh_{p-1}r^p/p$.
For any
$\bar{x}\in\R^p$ and any $\s>0$ we have
\begin{eqnarray*}
  \Prob_{x\sim N(\bar{x},\s^2\Id)}\big\{\|x\|\leq \e\big\}
  &\leq & \Prob_{x\sim N(0,\s^2\Id)}\big\{\|x\|\leq \e\big\}
  \;=\; \frac{1}{(\s\sqrt{2\pi})^p}
          \int_{\|x\|\leq \e} e^{-\frac{\|x\|^2}{2\s^2}} dx\\
  &\stackrel{\scriptstyle x=\s z}{=}&
    \frac{1}{(\sqrt{2\pi})^p}
          \int_{\|z\|\leq \frac{\e}{\s}}
         e^{-\frac{\|z\|^2}{2}} dz\\
  &\leq & \frac{1}{(\sqrt{2\pi})^p}\,
               \vol B_{\frac{\e}{\s}}
  \;=\; \frac{1}{(\sqrt{2\pi})^p}
               \Big(\frac{\e}{\s}\Big)^p\,\vol B_1\\
   &=& \frac{1}{(\sqrt{2\pi})^p}
               \Big(\frac{\e}{\s}\Big)^p\,
               \frac{\Oh_{p-1}}{p}.
\end{eqnarray*}
Taking $\bar{x}=\bar{a}_m^\perp$,
$\e=\frac{1}{\xi}$, and $p=n-m+1$
the claim~\eqref{eq:C4} follows.
\eproof

\proofof{Proposition~\ref{thm:main_tail}}
The proof is based on an idea in~\cite{sst:06}.
For $A\in\R^{m\times n}$ there exists $u_A\in\Sph^{m-1}$ such that
$\|A^{\dagger}\|=\|A^{\dagger}u_A\|$.
Moreover, for almost all $A$, the vector $u_A$ is uniquely determined
up to sign.
Using the singular value decomposition it is easy to show that,
for all $v\in\Sph^{m-1}$,
\begin{equation}\label{eq:sankar}
  \|A^{\dagger}v\| \geq \|A^{\dagger}\| \cdot |u_A\transp v|.
\end{equation}

Now take $A~\sim N(\oA,\s^2\Id)$ and
$v\sim U(\Sph^{m-1})$ independently.
Then, for any $s\in(0,1)$ and $t>0$ we have
\begin{align*}
  \Prob_{A,v}\big\{\|A^{\dagger}v\|\geq\, t & \sqrt{1-s^2}\big\}
  \;\geq\;
  \Prob_{A,v}\Big\{\|A^{\dagger}\|\geq t\ \&\
    |u_A\transp v|\geq \sqrt{1-s^2}\Big\}\\
  =\; & \Prob_{A}\big\{\|A^{\dagger}\|\geq t\big\} \cdot
     \Prob_{A,v}\Big\{|u_A\transp v|\geq \sqrt{1-s^2}
    \;\Big|\; \|A^{\dagger}\| \geq t\Big\}\\
  \geq\; & \Prob_{A}\big\{\|A^{\dagger}\|\geq t\big\} \cdot
   \sqrt{\frac{2}{\pi m}}\,s^{m-1},
\end{align*}
the last line by Lemma~\ref{lem:s} with $\xi=\sqrt{1-s^2}$.
Now we use Proposition~\ref{prop:bound1} with
$\xi=t\sqrt{1-s^2}$ to deduce that
\begin{eqnarray}\label{eq:bound2}
  \Prob_{A}\big\{\|A^{\dagger}\|\geq t \big\}
  &\leq& \sqrt{\frac{\pi m}{2}}\;
   \frac{1}{s^{m-1}}
  \Prob_{A,v}\{\|A^{\dagger}v\|\geq t\sqrt{1-s^2}\}\\
  &\leq&  \frac{\sqrt{m}}{2s^{m-1}}\;
       \frac{1}{(\sqrt{2\pi})^{n-m}}\,
               \frac{\Oh_{n-m}}{n-m+1}\,
               \Big(\frac{1}{\s t\sqrt{1-s^2}}\Big)^{n-m+1}.
               \nonumber
\end{eqnarray}

We next choose $s\in(0,1)$ to minimize the bound above.
To do so amounts to maximize
$(1-x)^{\frac{n-m+1}{2}} x^{\frac{m-1}{2}}$ where
$x=s^2\in(0,1)$, or yet, to maximize
$$
   g(x)=\Big((1-x)^{\frac{n-m+1}{2}}
    x^{\frac{m-1}{2}}\Big)^{\frac2n}
     =(1-x)^{\frac{n-m+1}{n}} x^{\frac{m-1}{n}}
     =(1-x)^{1-\lambda} x^\lambda .
$$
We have
$
   \frac{d}{dx}\ln g(x)= \frac{\lambda}{x}-\frac{1-\lambda}{1-x}
$
with the only zero attained at $x^*=\lambda$.

Replacing $s^2$ by $\lambda$ in~\eqref{eq:bound2} we obtain the bound
$$
  \Prob_{A}\big\{\|A^{\dagger}\|\geq t \big\} \ \leq\
  \frac{\sqrt{\lambda n+1}}
  {2\lambda^{\frac{\lambda n}{2}}}\;
       \frac{1}{(\sqrt{2\pi})^{n-m}}\,
               \frac{\Oh_{n-m}}{(1-\lambda)n}\,
               \bigg(\frac{1}{\s t\sqrt{1-\lambda}}\bigg)^{(1-\lambda)n} .
$$
Lemma~\ref{le:gbound} implies
$$
\lambda^{-\frac{\lambda n}{2}}
 = \Big( \lambda^{-\frac{\lambda}{2(1-\lambda)}}\Big)^{(1-\lambda)n}
 \le e^{\frac{(1-\lambda)n}{2}} .
$$
So we get
\begin{align*}
\Prob_{A}\big\{& \|A^{\dagger}\|\geq t\big\} \ \leq\
\frac{\sqrt{\lambda n+1}}{2}\;
       \frac{1}{(\sqrt{2\pi})^{n-m}}\,
               \frac{\Oh_{n-m}}{(1-\lambda)n}\,
               \bigg(\frac{\sqrt{e}}{\s t\sqrt{1-\lambda}}\bigg)^{(1-\lambda)n} \\
  \;=\;&
  \; \frac{\sqrt{\lambda n +1}}{2}\;
    \bigg(\frac{e}{1-\lambda}\bigg)^{\frac{(1-\lambda)n}{2}}\;
       \frac{1}{(\sqrt{2\pi})^{n-m}}\,
               \frac{\Oh_{n-m}}{(1-\lambda)n}\,
               \bigg(\frac{1}{\s t}\bigg)^{(1-\lambda)n}\\
 =\;& \frac1{2(1-\lambda)}\;
 \sqrt{\lambda +\frac1{n}}\;\frac1{\sqrt{n}}
\bigg(\frac{e}{1-\lambda}\bigg)^{\frac{(1-\lambda)n}{2}}\;
       \frac{\Oh_{n-m}}{(\sqrt{2\pi})^{n-m}}\,
               \bigg(\frac{1}{\s t}\bigg)^{(1-\lambda)n}\\
\le\;& \frac{\sqrt{\lambda + 1}}{2(1-\lambda)}\;\frac1{\sqrt{n}}
  \bigg(\frac{e}{1-\lambda}\bigg)^{\frac{(1-\lambda)n}{2}}\;
       \frac{2\pi^{\frac{n-m+1}{2}}}
       {\Gamma\big(\frac{n-m+1}{2}\big)(\sqrt{2\pi})^{n-m}}\,
               \bigg(\frac{1}{\s t}\bigg)^{(1-\lambda)n}\\
 = \;& \frac{\sqrt{1+\lambda}}{1-\lambda}\; \frac1{\sqrt{n}}\;
  \bigg(\frac{e}{1-\lambda}\bigg)^{\frac{(1-\lambda)n}{2}}\;
       \frac{\sqrt{2\pi}}
     {\Gamma\big(\frac{n(1-\lambda)}{2}\big)2^{\frac{(1-\lambda)n}{2}}}\,
               \bigg(\frac{1}{\s t}\bigg)^{(1-\lambda)n}.
\end{align*}
We next estimate $\Gamma\big(\frac{(1-\lambda)n}{2}\big)$. To do so,
recall Stirling's bound
$$
  \sqrt{2\pi}x^{x+\frac12}e^{-x}<\Gamma(x+1)
  <\sqrt{2\pi}x^{x+\frac12}e^{-x+\frac{1}{12x}}
 \qquad\mbox{for all $x>0$}
$$
which yields, using $\Gamma(x+1)=x\Gamma(x)$,
the bound $\Gamma(x) > \sqrt{2\pi/x}\, (x/e)^{x}$.
We use this with $x=\frac{(1-\lambda)n}{2}$ to obtain
$$
  \Gamma\Big(\frac{(1-\lambda)n}{2}\Big)
  \geq \sqrt{\frac{4\pi}{(1-\lambda)n}}\
  \Big(\frac{(1-\lambda)n}{2e}\Big)^\frac{(1-\lambda)n}{2}.
$$
Plugging this into the above we obtain
(observe the crucial cancellation of $\sqrt{n}$)
\begin{align*}
  \Prob_{A}\big\{&\|A^{\dagger}\|\geq t\big\} \\
  \leq\ & \sqrt{\frac{1+\lambda}{(1-\lambda)^2}}\; \frac1{\sqrt{n}}\;
  \bigg(\frac{e}{1-\lambda}\bigg)^{\frac{(1-\lambda)n}{2}}\;
     \, \sqrt{2\pi}\, \sqrt{\frac{(1-\lambda)n}{4\pi}}\,
   \Big(\frac{e}{(1-\lambda)n}\Big)^\frac{(1-\lambda)n}{2}
               \bigg(\frac{1}{\s t}\bigg)^{(1-\lambda)n}\\
  =\ &
  c(\lambda)\, \bigg(\frac{e}{1-\lambda}\bigg)^{(1-\lambda)n}\;
   \Big(\frac{1}{n}\Big)^{\frac{(1-\lambda)n}{2}}
               \bigg(\frac{1}{\s t}\bigg)^{(1-\lambda)n}
  = c(\lambda)\,\bigg(\frac{e}{\s\sqrt{n}(1-\lambda)t}\bigg)^{(1-\lambda)n}\;              ,
\end{align*}
which completes the proof of the proposition.
\eproof

\subsection{Proof of Theorem~\ref{thm:tailbound}}

To simplify notation we write
$c:=c(\lambda)$ and $Q:=Q(m,n)$.
Proposition~\ref{thm:main_tail} implies that for any $\e>0$ we have
\begin{equation}\label{eq:e_Adagger}
\Prob_{A\sim N(\oA,\s^2\Id)}\Big\{
  \|A^\dagger\|\geq \frac{e}{1-\lambda}\,
    \frac{1}{\s\sqrt{n}}\,
   \Big(\frac{c}{\e}\Big)^{\frac{1}{(1-\lambda)n}} \Big\}
   \ \leq\ \e.
\end{equation}
Similarly, letting $\e=e^{-\frac{t^2}{2\s^2}}$ in
Proposition~\ref{le:enorm}
and solving for $t$ we deduce that, for any $\e\in(0,1]$,
\begin{equation}\label{eq:e_A}
\Prob\Big\{\|A\|\geq  Q\s\sqrt{n}+\s\sqrt{2\ln\frac{1}{\e}} + 1
     \Big\}\ \leq\ \e.
\end{equation}
We conclude that
\begin{equation}\label{eq:tail_kappa}
   \Prob_{A\sim N(\oA,\s^2\Id)}\Big\{\kappa(A)\geq
   \frac{ez(\e)}{1-\lambda}\Big\}
   \ \leq\ 2\e,
\end{equation}
where we have have set, for $\e\in(0,1]$,
\begin{equation}\label{def:z(e)}
  z(\e):=
     \Bigg(Q + \sqrt{\frac{2}{n} \ln\frac1{\e}} + \frac1{\s\sqrt{n}}\Bigg)\,
        \Big(\frac{c}{\e}\Big)^{\frac{1}{(1-\lambda)n}}.
\end{equation}
We note that $z(1)=\z:=\z_\s(m,n)$, cf.~Equation~\eqref{eq:defB}.
Moreover, $\lim_{\e\to 0}z(\e) =\infty$ and $z$ is decreasing in the
interval $(0,1]$.
Hence, for $z\ge \z$, there exists $\e=\e(z)\in (0,1]$ such that $z=z(\e)$.

We need to upper bound $\e(z)$ as a function of $z$.
To do so, we start with a weak lower bound on $\e(z)$ and
claim that
\begin{equation}\label{eq:claim}
\frac1{n}\,\ln\frac1{\e}\ \le\ \ln (2z(\e)) .
\end{equation}
To show this, recall that
$Q\ge \sqrt{\frac{n}{n+1}} \ge \frac1{\sqrt{2}}$
due to Lemma~\ref{pro:Espnorm}.
Hence
$\z \ge Q \ge 1/\sqrt{2}$ and it follows that
$\sqrt{2} z \le 1$ for  $z \ge \z$.
Thus, Equation~\eqref{def:z(e)} implies that
$$
 z(\e)\ \ge\ \frac{1}{\sqrt{2}}\,
  \Big(\frac{c}{\e}\Big)^{\frac{1}{(1-\lambda)n}}.
$$
Using $c\ge \frac1{\sqrt{2}}$ we get
$$
 (\sqrt{2} z)^n \ \ge\ (\sqrt{2} z)^{(1-\lambda)n}
 \ge\ \frac{c}{\e}\ \ge\ \frac1{\sqrt{2}\,\e}.
$$
Hence $(2z)^n \ge 1/\e$, which shows
the claimed inequality~\eqref{eq:claim}.

Using the bound~\eqref{eq:claim} in Equation~\eqref{def:z(e)} we get,
again writing $z=z(\e)$, that
$$
 z\ \le\
 \Big(Q + \sqrt{2\ln (2z)} + \frac1{\s\sqrt{n}}\Big)\,
        \Big(\frac{c}{\e}\Big)^{\frac{1}{(1-\lambda)n}},
$$
which means
$$
 \e\ \le\
  c \bigg[
   \Big(Q + \sqrt{2\ln(2z)} + \frac1{\s\sqrt{n}}\Big)\,
    \frac1z\bigg]^{(1-\lambda)n}.
$$
By~\eqref{eq:tail_kappa} this completes the proof.
\eproof

\subsection{Proof of Corollary~\ref{cor:main}}

Fix $\lambda_0\in (0,1)$ and put $c:=c(\lambda_0)$.
Suppose that $m\le n$ satisfy $\lambda =(m-1)/n \le \lambda_0$.
Then $n-m +1 = (1-\lambda)n \ge (1-\lambda_0) n$ and
in order to have $n-m$ sufficiently large
it suffices to require that $n$ is sufficiently large.
Thus, $c^{\frac1{n-m+1}} \le 1.1$
if $n$ is sufficiently large. Similarly, because of
Lemma~\ref{pro:Espnorm},
$Q(m,n)\leq 2.1$
for large enough~$n$. This implies that,
for $\frac1{\sqrt{m}} \le \sigma\le 1$,
we have
$$
 Q(m,n)  + \frac1{\s\sqrt{n}}\ \le\
 2.1 + \frac1{\s\sqrt{n}} \ \le\
 2.1 + \sqrt{\frac{m}{n}} \ \le\
 2.1 + \sqrt{\lambda_0 +\frac1{n}}\ \le\ 3.1 ,
$$
provided $n$ is large enough. Then
$\z_\s(m,n) \le 3.1\cdot 1.1 = 3.41$.

By Theorem~\ref{thm:tailbound},
the random variable
$Z :=  (1-\lambda)\kappa(A)/e$
satisfies, for any $\oA$ with $\|\oA\|\le 1$ and
any $z\ge 3.41$,
\begin{eqnarray*}
 \Prob_{A\sim N(\oA,\s^2\Id)}\big\{ Z \ge z \big\}
  &\leq& 2c\,\bigg[\Big(Q(m,n)+ \sqrt{2\ln (2z)} + \frac1{\s\sqrt{n}}\Big)\,
    \frac1z\bigg]^{n-m+1} \\
  &\leq& 2c\,\bigg[\Big(3.1+ \sqrt{2\ln (2z)}\Big)\,
    \frac1z\bigg]^{n-m+1}.
\end{eqnarray*}
Since
$3.1 + \sqrt{2\ln (2z)}\le e\sqrt{z}$
for $z\ge 4$
we deduce that, for all such $z$,
$$
 \Prob_{A\sim N(\oA,\s^2\Id)}\big\{ Z \ge z \big\}
  \ \leq\ 2c\,\Big(\frac{e}{\sqrt{z}}\Big)^{n-m+1} .
$$
Using this tail bound to compute $\E(Z)$ we get
\begin{eqnarray*}
 \E(Z) &=& \int_0^\infty \Prob\{Z\ge z\}\, dz\
  \;\le\; e^2 + 2c
  \int_{e^2}^\infty
  \left(\frac{e^2}{z}\right)^{\frac{n-m+1}{2}}dz \\
  &\stackrel{\scriptstyle z=e^2y}{=}&
    e^2 + 2c\int_1^\infty
  \left(\frac{1}{y}\right)^{\frac{n-m+1}{2}}e^2dy
  \;=\; e^2+ \frac{4ce^2}{n-m-1}.
\end{eqnarray*}
We can now conclude since
$$
 \E((1-\lambda)\kappa(A))
  = \E(eZ) =e\E(Z)\leq e^3+ \frac{4ce^3}{n-m-1}
  \leq 20.1
$$
the inequality, again, by taking $n$ large enough.
\eproof

\section{Applications}
\label{sec:applications}

We next briefly discuss the two applications of our main result
mentioned in the introduction.

\subsection{Accuracy of Linear Least Squares}

Recall the problem~\eqref{eq:ls} described in the introduction, namely,
to compute the minimum of $\|Ax-b\|^2$ over $x\in\R^n$ 
for given $A\in\R^{m\times n}$ and $b\in\R^m$ (with $m>n$).
It is well known
that the loss of precision $\LoP(A^\dagger b)$
---that is, the number of correct digits in
the entries of the data $(A,b)$ minus the same number for the computed
solution $A^\dagger b$--- satisfies
(cf.~\cite{Wedin73} and \cite[Ch.~19]{higham:96}) 
$$
  \LoP(A^\dagger b)\leq \log m n^{3/2} + 2\log\kappa(A)+\Oh(1).
$$
 Corollary~\ref{cor:main}, combined with Jensen's inequality,
implies that $\E(\log\kappa(A)) \le \log(20.1/(1-\lambda)) = \Oh(1)$
under the assumptions stated in the corollary.
Hence for sufficiently elongated, large matrices $\oA$, the
expected loss of precision in the computation of the solution
$A^\dagger b$ over all small perturbations $A$ of $\oA$ is
dominated by the term $\log m n^{3/2}$.

\subsection{Complexity of the Conjugate Gradient Method}

If $P\in\R^{m\times m}$ is a symmetric positive definite
matrix and $c\in\R^m$, the system $Px=c$ can be solved by
the Conjugate Gradient Method (CGM), cf.~\cite{hest-stiefel:52}.
This is an iterative algorithm which performs at most $m$ iterations but
may require less. Indeed, it is known
(see, e.g., \cite[Lecture~38]{ThrefethenBau})
that an $\varepsilon$-approximation
of the solution~$x$ can be computed in at most
$\frac12\sqrt{\kappa(P)}|\ln\varepsilon| $
iterations
($\varepsilon$ measures the
relative error of the approximation with respect to the
Euclidean norm).

In many cases the matrix $P$ arises as $P=AA\transp$ for some
matrix $A\in\R^{m\times n}$ with $n>m$. 
If $A$ is standard Gaussian distributed, 
then the resulting distribution of $P$, 
called {\em Wishart distribution}, 
has been extensively studied in multivariate statistics. 
However, in our case of interest, $A$ is noncentered 
and much less is known about the resulting 
distribution of $P$ 
(called noncentral Wishart). 
Fortunately, using the fact that $\sqrt{\kappa(P)}=\kappa(A)$, 
we can directly apply our tail bounds for $\kappa(A)$ 
for a noncentral, isotropic Gaussian 
distribution of $A$,  
to derive bounds for the expected number
of iterations of CGM.

To do so we use again Corollary~\ref{cor:main}.
It shows that for all $\lambda_0\in(0, 1)$ and all
$0 <\sigma \le 1$ there exists $n_0$ such that for
all $1\leq m<n$ we have
$$
 \sup_{\|\oA\| \le 1} \E_{A\sim N(\oA,\sigma^2\Id)}(\kappa(A))
  \ \leq\ \frac{20.1}{1-\lambda},
$$
provided $\lambda = \frac{m-1}{n}\leq \lambda_0$ and $n\geq n_0$.
It follows that if $P$ is obtained as $AA\transp$ for a large, elongated,
rectangular matrix $A$ then, we should expect to
compute a solution with the desired accuracy with
about $\frac{1}{2}\frac{20.1}{1-\lambda}|\ln\varepsilon|$ iterations.
It is known that
each iteration of CGM takes $6n^2 +\Oh(n)$ arithmetic operations.
Therefore, the expected cost of running CGM on $P$ is
$$
  3n^2 \frac{20.1}{1-\lambda}|\ln\varepsilon|
 +\Oh(n)=\frac{60.3n^2}{1-\lambda}|\ln\varepsilon| +\Oh(n).
$$
The leading term in this expression is smaller than the
$\frac23 n^3$ operations performed
by Gaussian elimination as long as
$$
  \varepsilon\geq e^{-\frac{n(1-\lambda)}{91}}.
$$
For large $n$ (and $\lambda$ not too close to 1) this bound produces
very small values of $\varepsilon$
and therefore, CGM yields, on the average (both for a Wishart distribution
of data $P$ and for Wishart perturbations of arbitrary data), 
remarkably good approximations of the solution $x=P^{-1}c$.

\section{Some Numerical Simulations}
\label{sec:numerical}

Section~6 in~\cite{CDW:05} describes the result of numerical computations
producing experimental values for $\E(\ln\kappa(A))$ (for certain
choices of $\oA$ and $\sigma$), which are
denoted by $\Avr(\ln\kappa(A))$ and compared with
the upper bound for $\E(\ln\kappa(A))$
$$
  \bNA:=
  \ln\left(m
  +\sigma m\sqrt{5n}\right)
  +\ln \frac{2.35}{\sigma} + \frac{1}{r}
  +\sqrt{\frac{e\pi}{5}}
$$
obtained there. The data in
Tables~\ref{tab:2} to~\ref{tab:5}
is taken from~\cite{CDW:05}.
Each row in these tables corresponds to a pair $(m,n)$.
For each row, 500 random matrices $A\in\R^{m\times n}$
were computed following the distribution $N(\sqrt{m}\,\oA,\Id)$,
where $\oA$ was chosen as
$$
  \oA := \frac{{\sf ones(m,n)}}{\|{\sf ones(m,n)}\|},
$$
and {\sf ones(m,n)} denotes the $m\times n$ matrix all of whose
entries are 1.
The column with header $\bOA$ shows the empirical average of
$\ln\kappa(A)$ for the 500 chosen random matrices~$A$.
Since $\kappa(A)$ is scale invariant we note that this corresponds
to random matrices chosen from $N(\oA,\sigma\Id)$,
where $\sigma=1/\sqrt{m}$.

\begin{table}[h]
\begin{center}
\begin{tabular}{||c|c|c|c|c||}
 \hline
 $m$ & $n$ &$\bOA$ & $\bNA$ & $\ln (20.1/(1-\lambda))$\\
\hline \hline
10&15& 1.88278226808667 &7.73190477060415 & \\
\cline{1-4} 20&30&  2.04718612539162 & 8.74083698937094 &
\\
\cline{1-4}
40&60 &2.13539482051851&9.75820027818245& 4.0993321
\\
\cline{1-4}
80 & 120 & 2.19377719811291 & 10.78180469776403 &
\\
 \cline{1-4}
160& 240&2.23119383890675&11.80997066079053&
\\
\hline
\end{tabular}\caption{$n=1.5m.$}\label{tab:2}
\end{center}
\end{table}

\begin{table}[h]
\begin{center}
\begin{tabular}{||c|c|c|c|c||}
 \hline
 $m$ & $n$ &$\bOA$ & $\bNA$ & $\ln (20.1/(1-\lambda))$\\
\hline \hline
 5& 10 & 1.28204418194521&6.35902343647518&\\
\cline{1-4}
10&20& 1.48669849397793 &7.36178009761038&  \\
\cline{1-4}
20&40 &  1.59394635398509&8.37451330180407 & \\
 \cline{1-4}
 40& 80&1.64896402420115& 9.39470162365532& 3.693866\\
 \cline{1-4}
80&160&1.69565973841311&10.42037692088400&\\
 \cline{1-4}
 160 &320 &1.72154032592663 &11.45004561375610&\\
 \hline
\end{tabular} \caption{$n=2m.$}\label{tab:3}
\end{center}
\end{table}

\begin{table}[h]
\begin{center}
\begin{tabular}{||c|c|c|c|c||}
 \hline
 $m$ & $n$ &$\bOA$ & $\bNA$ & $\ln (20.1/(1-\lambda))$\\
\hline \hline
10&25& 1.24167342192086 &7.46370799208199&\\
\cline{1-4}
20&50 &1.34213347902230&8.47908853717777 &  \\
 \cline{1-4}
 40& 100&1.40120155287858&9.50123344342563& 3.511545 \\
 \cline{1-4}
80&200&1.44120596017225&10.52833707967242&\\
 \cline{1-4}
 160 &400 &1.45928497502137&11.55903912197539&\\
 \hline
\end{tabular} \caption{$n=2.5 m.$}\label{tab:4}
\end{center}
\end{table}

\begin{table}[h]
\begin{center}
\begin{tabular}{||c|c|c|c|c||}
 \hline
 $m$ & $n$ &$\bOA$ & $\bNA$ & $\ln (20.1/(1-\lambda))$\\  \hline \hline
 5& 15 &  0.98741849882614& 6.37209092337754 &  \\
\cline{1-4}
10&30&1.10550395287499 &7.38102314214432&  \\
\cline{1-4}
20&60 &1.18790345922560&8.39838643095583 & \\
 \cline{1-4}
 40& 120& 1.23914387557043&9.42199085053742& 3.406185 \\
 \cline{1-4}
80&240& 1.27096561714092&10.45015681356392&\\
 \cline{1-4}
 160 &480  &1.28600775609989&12.14829242876138&\\
 \hline
\end{tabular}\caption{$n=3m.$}\label{tab:5}
\end{center}
\end{table}

In~\cite{CDW:05} it is
observed that ``one sees that when one fixes $m$ and lets $n$ increase
the quantity $\Avr(\ln\kappa(A))$ decreases. This is in contrast
with the behaviour of $\mu(m,n,\sigma)$. It appears that our methods
are not sharp enough to capture the behaviour of
$\E(\ln\kappa(A))$.'' Compare now with the
results of the present paper. It follows from
Corollary~\ref{cor:main}, by Jensen's inequality, that,
for sufficiently large $n$,
$\E(\ln\kappa(A))\leq \ln \frac{20.1}{1-\lambda}$. But if $m$ is
held fixed then, when $n$ increases,
$\lambda$ decreases and so does $\ln \frac{20.1}{1-\lambda}$.

One still observes a difference between the bound $\ln \frac{20.1}{1-\lambda}$
and the values of $\bOA$. Part of this difference comes from the
asymptotic character of this bound and the fact that our data
is limited to $m\leq 160$. One sees on the tables that larger values
of $m$ would approach $\bOA$ to $\ln \frac{20.1}{1-\lambda}$.
We conjecture that, in addition to the possible loss of sharpness
coming from the use of Jensen's inequality, the difference above is due to
the roughness of the constant 20.1.

{\small

}
\end{document}